\documentclass[11pt, a4paper, oneside, reqno]{amsart}

\usepackage{bbm}

\newtheorem{theorem}{Theorem}[section]

\newtheorem{lemma}[theorem]{Lemma}

\newtheorem{remark}[theorem]{Remark}

\newtheorem{proposition}[theorem]{Proposition}

\DeclareMathOperator{\rank}{rank}

\makeatletter
\def\@settitle{\begin{center}%
		\baselineskip14\p@\relax
		\normalfont\LARGE\scshape\bfseries
		\@title
	\end{center}%
}
\makeatother

\makeatletter

\def\subsection{\@startsection{subsection}{2}%
	\z@{.5\linespacing\@plus.7\linespacing}{.5\linespacing}%
	{\normalfont\large\bfseries}}

\def\subsubsection{\@startsection{subsubsection}{3}%
	\z@{.5\linespacing\@plus.7\linespacing}{.5\linespacing}%
	{\normalfont\itshape}}

\usepackage{caption,subcaption}
\usepackage[colorlinks	= true,
			raiselinks	= true,
			linkcolor	= MidnightBlue,
			citecolor	= Mahogany,
			urlcolor	= ForestGreen,
			pdfauthor	= {Peyman Mohajerin Esfahani},
			pdftitle	= {},
			pdfkeywords	= {},   
			pdfsubject	= {},
			plainpages	= false]{hyperref}
\usepackage[usenames, dvipsnames]{color}
\usepackage{dsfont,amssymb,amsmath, graphicx,fancyhdr,mdframed}
\usepackage{amsfonts,dsfont,mathtools, mathrsfs,amsthm,wrapfig} 
\usepackage[]{siunitx}
\usepackage{ragged2e}
\usepackage{enumitem}

\allowdisplaybreaks
\usepackage{algorithm}
\usepackage[noend]{algpseudocode}
\usepackage[square,numbers]{natbib}

\allowdisplaybreaks
\date{\today}

\patchcmd\@setauthors
  {\MakeUppercase{\authors}}
  {\authors}
  {}{}

\graphicspath{{../fig/}}

\newcommand{\Ball}[1]{\mathbb{B}_{\varepsilon}^{#1}(\widehat{\Sigma})}
\DeclareMathOperator{\Trace}{Tr}
\newcommand{\Tr}[1]{\Trace\left(#1\right)}
\newcommand{\Avg}[1]{\overline{#1}}

\DeclareMathOperator{\dist}{d}
\DeclareMathOperator{\Fro}{F}
\DeclareMathOperator{\KL}{KL}
\DeclareMathOperator{\Gel}{G}

\newcommand{\Lipschitz}{\mathcal{L}}

\newcommand{\SPlusSet}{\mathbb{S}_{+}}
\newcommand{\DPlusSet}{\mathbb{D}_{+}}

\newcommand{\SSetProj}[1]{\mathbb{S}_{#1}}

\newcommand{\dualFunc}[1]{g(#1)}

\newcommand{\LMO}{O}

\title[A Saddle Point Algorithm for Robust Data-Driven Factor Model Problems]{A Saddle Point Algorithm for \\ Robust Data-Driven Factor Model Problems}

\author[S. Khodakaramzadeh]{Shabnam Khodakaramzadeh$^1$}
\author[S. {Shafiee}]{Soroosh {Shafiee$^2$}}
\author[Gabriel de A.~Gleizer]{Gabriel {de Albuquerque Gleizer}$^1$}
\author[P. {Mohajerin Esfahani}]{Peyman {Mohajerin Esfahani$^{1,3}$}}

\thanks{The authors are with (1) Delft University of Technology, (2) Cornell University, and (3) the University of Toronto. This research is supported by the ERC Starting Grant TRUST-949796 and by the NSF CAREER ECCS-2541066.}

\begin{document}

\maketitle

\begin{abstract} 
We study the factor model problem, which aims to uncover low-dimensional structures in high-dimensional datasets. Adopting a robust data-driven approach, we formulate the problem as a saddle-point optimization. Our primary contribution is a first-order algorithm that solves this reformulation by leveraging a linear minimization oracle (LMO). We further develop semi-closed form solutions (up to a scalar) for three specific LMOs, corresponding to the Frobenius norm, Kullback-Leibler divergence, and Gelbrich (aka Wasserstein) distance. The analysis includes explicit quantification of these LMOs' regularity conditions, notably the Lipschitz constants of the dual function, which govern the algorithm's convergence performance. Numerical experiments confirm our method's effectiveness in high-dimensional settings, outperforming standard off-the-shelf optimization solvers.
\end{abstract}


\section{Introduction}\label{sec:introduction}
The correlation structure among random variables can be expressed using factor analysis in terms of common \emph{factors}~\cite{ComputationTheMaxLikelihood_khamaru}, i.e., the information in high-dimensional data can be compressed into unobserved factors~\cite{FactorAnalysisofMovingAverage_Zorzi}. Mathematically, a high-dimensional point $\xi \in \mathbb{R}^{n}$ is represented as sum of two independent unobserved parts
\begin{equation}\label{Eq_LinRegVec} 
    \xi = \Phi \alpha + \omega,
\end{equation}
where {\em factor loading matrix} $\Phi\in \mathbb{R}^{n \times r}$ is tall~($n \gg r$) and full-rank. The vector $\alpha \in \mathbb R^r$, with covariance $\Sigma_\alpha = I$, contains independent latent factors. Hence, $\Phi \alpha$, having interrelated components, specifies the low-dimensional representation of $\xi$. The \emph{idiosyncratic noise} $\omega$ has independent components~\cite{Multivar_Statistic_Analysis_Anderson,Factor_Models_Zorzi}. Such low-dimensional representations have been also studied in principal component analysis (PCA), which is well-suited for scenarios where the data is corrupted by small, unstructured noise~\cite{Pattern_Recognition_Bishop}, and also compressed sensing, where the noise is small and unstructured, but $r$ factors are contained in a sparse tall vector~$\alpha$~\cite{donoho2006compressed}. In contrast, the noise in the factor model has independent components and may be substantial. Originating from the analysis of mental test scores~\cite{Multivar_Statistic_Analysis_Anderson,GeneralIntelligence_Spearman}, the factor model has found broad applications across various domains, including control and system identification~\cite{ModelingComplexSystemsbyGFA_Picci,Hidden_Factor_Estimation_Zorzi,A_Robust_Approach_to_ARMA_FM_Zorzi}, fault and anomaly detection~\cite{MultirateFactorAnalysisModelsForFaultDetection_Zhou,FactorAnalysisBasedAnomaly_Wu}, econometrics~\cite{TheGeneralizedDynamicFM_Forni_2,AOneFactorMultivariateTimeSeries_Watson,StochasticRealizationProblems_Schuppen}, and statistics~\cite{Factor_Models_Zorzi}. When the covariance of $\xi$, denoted by $\Sigma$, is available, and assuming both random variables~$\xi$ and $\omega$ are zero-mean and independent from each other, we can rewrite~\eqref{Eq_LinRegVec} as
\begin{equation} \label{Eq_Cov}
    \Sigma = L + D,
\end{equation}
where $L = \Phi \Sigma_\alpha \Phi^\top$ is low rank with $\rank(L) = r$~(i.e., number of factors) and the non-negative diagonal matrix $D$ is the noise covariance. Often, $\Sigma$ is only observable by a finite dataset $\{\xi_k\}_{k=1}^N$, and is thus approximated with
\begin{equation}\label{Eq_SigmaHat}
    \widehat{\Sigma}= \frac{1}{N} \sum_{k=1}^N \big(\xi_k - \widehat{\mu}\big) \big(\xi_k - \widehat{\mu}\big)^\top, \quad \widehat{\mu} =  \frac{1}{N} \sum_{k=1}^N \xi_k.  
\end{equation}	
One can consider $\widehat{\mu} = 0$, under the assumption that~$\xi$ is zero-mean. To robustify to this approximation error, a family of covariance matrices is considered around $\widehat{\Sigma}$ as
\begin{align}\label{ball}
	 \Ball{\dist} \coloneqq \big \{ \Sigma \succeq 0 ~ \colon ~ \dist(\Sigma,\widehat{\Sigma}) \leq \varepsilon \big \},
\end{align}
 where $\dist$ is a generic distance function in the space of matrices, and $\varepsilon$ is the radius of the set. Therefore, our {\em robust data-driven factor model} problem is formulated as
 \begin{align} \label{Eq_FM_Convex_App} 
	\begin{array}{rcl} 
		J^\star  \coloneqq & \min\limits_{L, D } & \Tr{L} \\
		&\textrm{s.t.} & L \in \SPlusSet, D \in \DPlusSet,  L + D \in  \Ball{\dist}
	\end{array}
\end{align}
where $\SPlusSet$ is the cone of positive semidefinite (PSD) matrices, $\DPlusSet$ is the cone of non-negative diagonal matrices. The last constraint is equivalent to $\dist(L + D,\widehat{\Sigma}) \leq \varepsilon$. The trace operator is a standard convexification for the rank function~\cite{recht2010guaranteed}, i.e., the objective is to find the \emph{least number} of factors explaining the data. The hyperparameter~$\varepsilon$ depends on the precision of the approximation~$\widehat{\Sigma}$.

Traditionally, the factor model problem was addressed assuming that $\widehat{\Sigma}$ is an accurate estimate of $\Sigma$ (i.e., $\varepsilon = 0$)~\cite{mclachlan2007algorithm,Pattern_Recognition_Bishop}. However, recent research accounts for uncertainty in $\widehat{\Sigma}$, (i.e., $\varepsilon > 0$). For example, \cite{Alternating_Min_Zorzi} proposes a coordinate descent-type algorithm to minimize $\| \widehat{\Sigma} - L - D \|_{\Fro}^2$, while \cite{Certifiably_Optimal_Low_Rank_Bertsimas} develops an algorithm based on conditional gradient method targeting $\| \widehat{\Sigma} - L - D \|_{\mathrm{q}}^q$. The robust factor model problem~\eqref{Eq_FM_Convex_App} and its dynamic counterpart, limited to special cases where $\dist$ is the Kullback-Leibler (KL) and Itakura–Saito
divergences, are studied in~\cite{Factor_Models_Zorzi} and \cite{A_Robust_Approach_to_ARMA_FM_Zorzi}, respectively. However, unlike \cite{Alternating_Min_Zorzi,Certifiably_Optimal_Low_Rank_Bertsimas,Factor_Models_Zorzi,A_Robust_Approach_to_ARMA_FM_Zorzi}, which are limited to a specific choice of~$\dist$, we introduce an algorithm for the generic~$\dist$, which only requires access to the linear minimization oracle (LMO)
\begin{align} \label{Eq_LMO}
    \LMO(\Lambda) \coloneqq \arg\min_{\Sigma} \left\{ \langle \Lambda , \Sigma \rangle: \Sigma \in \Ball{\dist} \right\}
\end{align}
for any symmetric matrix $\Lambda$\footnote{If the program~\eqref{Eq_LMO} has multiple optimizers, the oracle~$\LMO(\Lambda)$ can be an arbitrary selection from the set.}. From a computational perspective, this study focuses on developing an efficient algorithm tailored to the semidefinite programming (SDP)~\eqref{Eq_FM_Convex_App}, for which the commercial solvers (e.g., MOSEK~\cite{SemidefOpt_MOSEK}) typically rely on second-order methods, and hence do not scale well to large problems. Moreover, General-purpose first-order methods (e.g., SCS~\cite{odonoghue2016conic}) instead require projection oracles that involve solving quadratic objectives over the feasible set, in contrast to the LMO~\eqref{Eq_LMO}, which only requires a linear objective. In fact, the LMO \eqref{Eq_LMO} can admit closed-form solutions, thereby avoiding the need to solve a full SDP. It is also worth noting that the problem~\eqref{Eq_FM_Convex_App} involves conic constraints associated with $\SPlusSet$ and $\DPlusSet$, which can make the direct solution computationally demanding.

In this work, we consider the static factor model problem, where the factors remain constant over time; for a dynamic counterpart,  where the factors temporal evolution is modeled via state-space equations, see, e.g.,~\cite{Hidden_Factor_Estimation_Zorzi}.
   
We summarize the contributions of this work as follows.
\begin{enumerate}[label=(\roman*), itemsep = 1mm, topsep = 0mm, leftmargin = 4mm]
    \item {\bf Saddle point characterization.}
    For possible conic structural information about the covariance matrix and generic distance functions, we reformulate factor model~\eqref{Eq_FM_Convex_App} as a saddle point problem~(Proposition~\ref{Prop_saddle}).
    
    \item {\bf First-order algorithm \& convergence.}
    Leveraging the saddle point reformulation and given an LMO, we propose a first-order algorithm with convergence guarantees derived based on standard regularity conditions such as Lipschitz constants~(Proposition~\ref{Prop_Alg}). A particular feature of the proposed algorithm is the linear convergence rate of the projection operator, enabled by Dykstra's projection technique, as opposed to the standard sublinear rate~(Proposition~\ref{prop_projection}).

    \item {\bf Special LMOs: closed-form description \& Lipschitz constant.} The choice of distance function influences the LMO and Lipschitz constant of the dual function, both of which play a central role in the proposed algorithm convergence. We derive the closed-form description for the LMO and its respective Lipschitz constant for the special cases of Frobenius norm~(Proposition~\ref{Prop_Oracle_Lipschitz_FroNormConst}), KL divergence~(Proposition~\ref{Prop_Oracle_Lipschitz_KLConst}), and Gelbrich distance~(Proposition~\ref{Prop_Gelbrich}). As a special byproduct, we show that the Gelbrich distance is strongly convex with respect to the Frobenius norm, a property of interest for optimization algorithms~(Remark~\ref{rem:strong convexity}). In comparison with the existing literature, to our best knowledge, this appears to be the first result that also applies to low-rank matrices.
\end{enumerate}
The theoretical results are validated through extensive numerical experiments, demonstrating the performance of the proposed algorithm. To improve the flow and readability,  all the technical proofs are relegated to Appendix~\ref{Append_Proofs}. To facilitate reproducibility, we provide an open-source MATLAB library available at \url{https://github.com/skhodakaram/Factor_Model}.

\textbf{Notation.} For any symmetric matrix $A$, its maximum eigenvalue is represented as $\lambda_{\textrm{max}} (A)$, its trace is denoted by $\Tr{A}$, and the vector $\lambda^{\frac{1}{2}}(A)$ contains the square roots of all its eigenvalues. The diagonal matrix $\textrm{Diag}(A)$ contains the diagonal entries of $A$. The inner product of any $A, B \in \mathbb{R}^{n \times m}$ is denoted as $\langle A , B \rangle \coloneqq \Tr{A^\top B}$. The Frobenius and nuclear norms of $A\in \mathbb{R}^{n \times m}$ are denoted as $\|A\|_{\Fro} \coloneqq \sqrt{\langle A,A \rangle}$ and $\| A \|_{\ast}$, respectively. The Euclidean norm of $a \in {\mathbb R}^n$ is denoted by~$\|a\|_2$. The relative interior of a set $\mathbb{A}$ is represented as ${\text{\normalfont rint}}(\mathbb{A})$. The operator~$\Pi_{\mathbb{A}} \left[x \right] := \arg\min_{y \in \mathbb{A}} \| x - y \|_F^2$ denotes the orthogonal projection of the point $x$ onto the set $\mathbb{A}$. The element-wise inequality between matrices is represented as $A \geq B$, and the semidefinite counterpart as $A \succeq B$. The positive semidefinite~(PSD) cone is denoted by~$\SPlusSet$. The dual cone of a convex cone $\mathbb{C} \in \mathbb{R}^n$ is represented as $\mathbb{C}^\ast \coloneqq \{x \in \mathbb{R}^n:\langle x, y \rangle \ge 0, \forall y \in \mathbb{C} \}$. The space of all lower triangular matrices in $\mathbb{R}^{n \times n}$ is denoted by $\mathbb{L}_n$. The normal cone of a set~$\mathbb{A}$ at~$x$ is defined by~$\mathcal{N}_{\mathbb{A}}(x) := \{ v : \langle v, y - x \rangle \leq 0, ~~ \forall y \in \mathbb{A} \}$. The indicator function over the set $\mathbb{A}$ is denoted as~$\mathbbm{1}_{\mathbb{A}}$.    

\section{Saddle Point Reformulation and First-Order Algorithm} \label{Sec_Primal_Dual_Alg}
   
\subsection{Saddle point characterization}
The first result of this paper is a saddle point (max-min) reformulation of the factor model problem~\eqref{Eq_FM_Convex_App}. This reformulation paves the way for optimization algorithms, especially considering the availability of the LMO \eqref{Eq_LMO}. 

\begin{proposition}[Saddle point reformulation]\label{Prop_saddle}
The optimal value $J^\star$ of the factor model problem in \eqref{Eq_FM_Convex_App} is equivalent to the max-min problem
\begin{align} \label{Eq_Our_FM_GeneralDist} 
    J^\star= \max\limits_{\renewcommand{\arraystretch}{1.0} \begin{matrix}
                         \scriptstyle I-\Lambda \in \SPlusSet \\
                           \scriptstyle  -\Lambda \in \DPlusSet^\ast
                          \end{matrix}}
            \min\limits_{\renewcommand{\arraystretch}{1.0}  \begin{matrix}
                        \scriptstyle \Sigma \in \Ball{\dist}
                          \end{matrix}} 
                          \langle \Lambda , \Sigma \rangle,  
\end{align}
where $\DPlusSet^\ast = \big\{A : {\rm Diag}(A) \ge 0\big\}$ is the dual cone of~$\DPlusSet$. 
\end{proposition}

We note that LMO \eqref{Eq_LMO} is used to find the solution to the inner minimization of \eqref{Eq_Our_FM_GeneralDist} as a function of the decision variable of the outer maximization problem. The max-min problem~\eqref{Eq_Our_FM_GeneralDist} facilitates development of an optimization algorithm to tackle \eqref{Eq_FM_Convex_App} numerically. The inner minimal value, hereafter referred to as the dual function, is 
\begin{align} \label{Eq_dLambda}
    \dualFunc{\Lambda} \coloneqq \min_{\Sigma \in \Ball{\dist}} \langle \Lambda, \Sigma \rangle 
\end{align}
Note that $\dualFunc{\Lambda}$ in \eqref{Eq_dLambda} is indeed the optimal value corresponding to the optimal point in LMO \eqref{Eq_LMO}. The availability of this oracle motivates us to develop a first-order algorithm for \eqref{Eq_dLambda}. The Lipschitz continuity of $\dualFunc{\Lambda}$ in \eqref{Eq_dLambda} is a critical regularity condition to ensure the success of the algorithm. For instance, if $\dualFunc{\Lambda}$ is Lipschitz continuous and its (sub)gradient is available, then one can use the classical projected gradient ascent with the stepsize proportion to $1/\sqrt{t}$, where $t$ is the iteration count~\cite{shamir2013stochastic}. The next result quantifies the Lipschitz constant of $\dualFunc{\Lambda}$. 

\begin{lemma} [Dual function Lipschitz constant] \label{lemma_LipschitzContinuity_dLambda}
The function $\dualFunc{\Lambda}$, defined in \eqref{Eq_dLambda}, is Lipschitz continuous with the constant $\Lipschitz$, i.e.,
\begin{align} \label{Eq_LipschitzConstant_general}
    \begin{array}{c}
    |\dualFunc{\Lambda_1} - \dualFunc{\Lambda_2}| \leq \Lipschitz \|\Lambda_1 - \Lambda_2\|_{\Fro}, \quad
    \text{where} \quad \Lipschitz \coloneqq \max\limits_{\Sigma \in \Ball{\dist}} \| \Sigma \|_{\Fro}.
    \end{array} 
\end{align} 
\end{lemma}

\subsection{First-order algorithm and Dykstra projection oracle} \label{Sec_FirstOrderAlg}
 We propose a first-order algorithm using the LMO \eqref{Eq_LMO} to solve the saddle point problem \eqref{Eq_Our_FM_GeneralDist}. 

\begin{proposition}[Algorithm \& convergence]\label{Prop_Alg}
Consider the optimization algorithm
\begin{equation} \label{Eq_OurFirstOrderAlgorithm}
    \left\{\begin{array}{rcl}
      \Sigma_t &\!\! = & \LMO(\Lambda_t) \\
      \Lambda_{t+1} & = & \Pi_{\SSetProj{1} \cap \SSetProj{2}} \left[ \Lambda_t + 
    \delta \Sigma_t \right] \\
      \Avg{\Lambda}_t & = & \frac{t-1}{t} \Avg{\Lambda}_{t-1} + \frac{1}{t} \Lambda_t, \quad \Avg{\Sigma}_t = \LMO(\Avg{\Lambda}_t)
    \end{array}
    \right. 
\end{equation}
where $\SSetProj{1} = \{ \Lambda : {\rm Diag}(\Lambda) \le 0 \}$, $ \SSetProj{2} = \{ \Lambda : I - \Lambda \in \SPlusSet \}$ are the conic constraint sets of $\Lambda$ in~\eqref{Eq_Our_FM_GeneralDist}, $\delta$ is a constant stepsize, and $\LMO$ is the LMO~\eqref{Eq_LMO}. Starting from any symmetric initial condition~$\Lambda_1$, after $T$ iterations we have
    \begin{align}\label{Alg_bound}
        0\le \langle \Lambda^\star, \Sigma^\star \rangle - \langle \Avg{\Lambda}_T, \Avg{\Sigma}_T \rangle \leq \frac{\| \Lambda_1 - \Lambda^\star \|^2}{2 \delta T} + {\frac{\delta}{2}} \Lipschitz^2
    \end{align}
where $(\Lambda^\star, \Sigma^\star)$ is a saddle point solution of~\eqref{Eq_Our_FM_GeneralDist}, and the constant~$\Lipschitz$ is the Lipschitz constant defined in~\eqref{Eq_LipschitzConstant_general}.  
\end{proposition}

\begin{remark}[Factor decomposition of $\Avg{\Sigma}_T$]
    The solution~$\Avg{\Sigma}_T$, obtained by \eqref{Eq_OurFirstOrderAlgorithm}, can be decomposed to $L+D$ using the standard factor analysis techniques from~\cite{mclachlan2007algorithm,Pattern_Recognition_Bishop}, which is substantially computationally cheaper than its robust counterpart \eqref{Eq_FM_Convex_App}. Let us recall that the matrix~$D$ gives information about the idiosyncratic noise variances, and the Cholesky factorization of the low-rank matrix $L$ provides the factor loading matrix $\Phi$.
\end{remark}

The next remark explains the choice of stepsize in \eqref{Eq_OurFirstOrderAlgorithm}.

\begin{remark}[Diminishing stepsize \& averaging]
    Minimizing the error bound in \eqref{Alg_bound} using the constant stepsize~$\delta$ reveals that the optimal choice is $\delta = \mathcal{O}(1 / \sqrt{T})$, yielding a suboptimality gap of $\mathcal{O}(1 / \sqrt{T})$; see \cite[Corollary~2]{bach2014adaptivity}. However, constant stepsizes often perform poorly in practice. Hence, we adopt a diminishing \emph{anytime} stepsize $\delta_t = \mathcal{O}(1/\sqrt{t}),$ which is widely used and achieves similar convergence, up to a logarithmic factor \cite{shamir2013stochastic}. 
\end{remark}

A possible difficulty in Algorithm~\eqref{Eq_OurFirstOrderAlgorithm} is the projection onto~$\SSetProj{1} \cap \SSetProj{2}$ for a given symmetric~$\bar \Lambda$. To address this issue, we leverage Dykstra's algorithm defined as follows: 
\vspace{0.5mm}
\begin{algorithm}[H]
    \caption{Dykstra's projection ($\Pi_{\SSetProj{1} \cap \SSetProj{2}}[\Lambda]$)}\label{Alg_DykstraProjection}
      \begin{algorithmic}[1]
        \State \textbf{Input:} $U_2^0= \Lambda$, $Z_1^0= 0$, $Z_2^0= 0$, threshold $\epsilon$
        \While{$\frac{\|U_2^k - U_1^k\|_{\Fro}}{\|U_2^k\|_{\Fro}} \geq \epsilon$} 
          \State $U_1^k= \Pi_{\SSetProj{1}} [U_2^{k-1} + Z_1^{k-1}]$ 
          \State $Z_1^k= U_2^{k-1} + Z_1^{k-1} - U_1^k$ 
          \State $U_2^k= \Pi_{\SSetProj{2}} [U_1^k + Z_2^{k-1}]$ 
          \State $Z_2^k= U_1^k + Z_2^{k-1} - U_2^k$
          \State $k \gets k+1$
        \EndWhile
        \algorithmicreturn \hspace{0.1mm} $U_2^k$
      \end{algorithmic}
\end{algorithm}

We know that Algorithm~\ref{Alg_DykstraProjection} converges \emph{asymptotically} under mild assumptions~\cite[Theorem~2]{boyle1986method}. However, this result can be improved to an exponential rate (aka {\em linear convergence}) for a certain class of conic problems if a corresponding dual admits the Kurdyka-{\L}ojasiewicz property~\cite[Theorem~5.3]{wang2024convergence}. The next result shows that these conditions are indeed satisfied in our setting.

\begin{figure}
    \centering
    \includegraphics[width=0.4\linewidth]{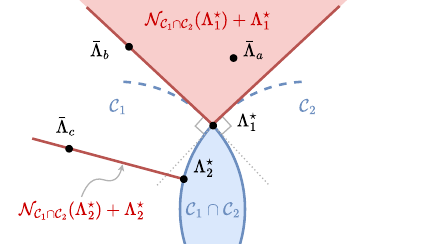}
    \caption{Depiction of the assumption in Proposition~\ref{prop_projection}. Any point in the relative interior of the normal cone (in red), such as $\bar\Lambda_a,$ satisfies the condition, unlike $\bar{\Lambda}_b$ and~$\bar{\Lambda}_c$.}
    \label{fig:normal cone}
\end{figure}

\begin{proposition}[Linear convergence for projection]
\label{prop_projection}
Let $\Lambda^\star = \Pi_{\SSetProj{1} \cap \SSetProj{2}}(\bar\Lambda)$, and assume that $\bar\Lambda - \Lambda^\star$ resides in the relative interior of the normal cone of $\SSetProj{1} \cap \SSetProj{2}$, i.e., $\bar\Lambda - \Lambda^\star \in {\text{\normalfont rint }}(\mathcal{N}_{\SSetProj{1} \cap \SSetProj{2}}(\Lambda^\star))$ (cf. Figure~\ref{fig:normal cone}). Then, Algorithm~\ref{Alg_DykstraProjection} converges linearly to $\Lambda^\star$ from any initial point.
\end{proposition}

\begin{remark} [Relative interior]
Condition~$\bar\Lambda - \Lambda^\star \in {\text{\normalfont rint}}(\mathcal{N}_{\SSetProj{1} \cap \SSetProj{2}}(\Lambda^*))$ can be viewed as a genericity condition on $\bar\Lambda$. Indeed, we have $\bar\Lambda- \Lambda^\star \in \mathcal{N}_{\SSetProj{1} \cap \SSetProj{2}}(\Lambda^\star)$ for any $\bar\Lambda$~\cite[Prop.~6.47]{bauschke2017convex}. Thus, the assumption excludes that $\bar \Lambda - \Lambda^\star$ belongs to the boundary, which almost never occurs when $\mathcal{N}_{\SSetProj{1} \cap \SSetProj{2}}(\Lambda^\star)$ is not a ray; see Figure~\ref{fig:normal cone}. If it is a ray, such as $\bar{\Lambda}_c$ in Figure~\ref{fig:normal cone}, then the projection to one of the sets solves the problem, which enables avoiding Dykstra's projection. 
\end{remark}
  
\section{Special Cases of Linear Minimization Oracle} \label{Sec_LMO_SpecialCases} 
    In this section, we study three special cases, the Frobenius distance, the KL divergence, and the Gelbrich distance, of LMO \eqref{Eq_LMO}, in particular in view of their computational complexity and Lipschitz continuity of~\eqref{Eq_dLambda}.
 
\subsection{Frobenius norm} 

The first distance is the Frobenius norm $\Fro(\Sigma,\widehat{\Sigma}) \coloneqq \| \Sigma-\widehat{\Sigma} \|_{\Fro}$. We show that in this case, the LMO~\eqref{Eq_LMO} admits an explicit form up to a scalar convex optimization.

\begin{proposition}[Frobenius oracle \& Lipschitz constant]\label{Prop_Oracle_Lipschitz_FroNormConst}
Consider LMO \eqref{Eq_LMO} for a given $\widehat{\Sigma}$ where the distance function is the Frobenius norm $\dist(\Sigma,\widehat{\Sigma}) = \Fro(\Sigma,\widehat{\Sigma})$.
\begin{enumerate} [label=(\roman*), leftmargin = 4mm]
\item\label{Prop_Oracle_Fro}\textbf{Closed-form description:} For any matrix $\Lambda$ and a positive scalar $\gamma$, we define
  \begin{equation}\label{Eq_LMO_Fro_PrimalDual_SigmaS}
    \Sigma^\star(\Lambda,\gamma) \coloneqq \Pi_{ \succeq 0 } \left[ \widehat{\Sigma} - \frac{1}{2 \gamma} \Lambda \right], 
  \end{equation}  
where $\Pi_{\succeq 0}$ denotes the projection with respect to the Frobenius norm onto the PSD cone. Then, $\LMO(\Lambda)$ in~\eqref{Eq_LMO} equates to $\Sigma^\star(\Lambda,\gamma^\star)$, where $\gamma^\star$ is the scalar solution to 
\begin{align} \label{Eq_LMO_Fro_gamma}
    \max_{0 < \gamma \leq \frac{1}{\varepsilon} \| \Lambda \|_{\Fro}} \langle \Lambda , \Sigma^\star(\Lambda,\gamma) \rangle + \gamma \big(\| \Sigma^\star(\Lambda,\gamma) - \widehat{\Sigma} \|_{\Fro}^2 - \varepsilon^2 \big). 
\end{align} 
\item\label{Prop_Lipschitz_Fro} \textbf{Lipschitz constant:} The Lipschitz constant of $\dualFunc{\Lambda}$ in Lemma \ref{lemma_LipschitzContinuity_dLambda} is bounded by
\begin{equation} \label{Eq_LipschitzConstantFro}
\Lipschitz \le \varepsilon + \| \widehat{\Sigma} \|_{\Fro}
\end{equation}
\end{enumerate}   
\end{proposition}
          
Proposition~\ref{Prop_Oracle_Lipschitz_FroNormConst}\ref{Prop_Oracle_Fro} provides an efficient computational way to implement the LMO \eqref{Eq_LMO} since the scalar concave maximization problem~\eqref{Eq_LMO_Fro_gamma} can be efficiently solved using bisection within the feasible region $\gamma \in (0,\frac{1}{\varepsilon} \| \Lambda\|_{\Fro}]$.
          
          \subsection{Kullback-Leibler divergence}
     
          The second case is the KL divergence between two identical mean normal distributions with covariance matrices $\Sigma$ and $\widehat{\Sigma}$~\cite[Chap. 2]{ref:cover1999elements}, defined as
          \begin{equation}\label{Eq_KL}
              \KL(\Sigma || \widehat{\Sigma}) \coloneqq \frac{1}{2} \Big( - \log \det{\Sigma} + \log \det{\widehat{\Sigma}} + \Tr{\Sigma \widehat{\Sigma}^{-1}} - n \Big),
          \end{equation}
          To make \eqref{Eq_KL} well defined, we assume $\widehat{\Sigma} \succ 0$, and thus invertible. We also assume $\xi$, $\alpha$, and $\omega$ are Gaussian random vectors. The result can be extended to the case of rank deficient $\widehat{\Sigma}$, where the KL-ball is required to have the same null space. While the KL divergence is a similarity measure, it is not symmetric in its arguments and does not satisfy the triangle inequality. Next, we provide an explicit description of its oracle and dual function Lipschitz constant.

\begin{proposition} [KL oracle \& Lipschitz constant] \label{Prop_Oracle_Lipschitz_KLConst}
Consider LMO \eqref{Eq_LMO} for a given $\widehat{\Sigma}$ where the distance function is the KL divergence $\dist(\Sigma,\widehat{\Sigma}) = \KL(\Sigma||\widehat{\Sigma})$.
    \begin{enumerate} [label=(\roman*), leftmargin = 4mm]
    \item\textbf{Closed-form description:} \label{Prop_Oracle_KL}
    For any matrix $\Lambda$ and a positive scalar $\gamma$, we define
    \begin{subequations}
    \begin{equation} \label{Eq_KL_SigmaStar}
        \Sigma^\star(\Lambda,\gamma) \coloneqq \left( \widehat{\Sigma}^{-1} + \frac{2}{\gamma} \Lambda \right)^{-1}.
    \end{equation}
    Then, $\LMO(\Lambda)$ in \eqref{Eq_LMO} equates to $\Sigma^\star(\Lambda,\gamma^\star)$, where the scalar $\gamma^\star$ satisfies the equations
    \begin{align} \label{Eq_KL_gamma_update}
    \begin{cases}
     &\KL ( \Sigma^\star (\Lambda,\gamma^\star) || \widehat{\Sigma} ) - \varepsilon = 0 \\
     &\begin{array}{l}
     \max \Big \{0 , 2 \lambda_{\textrm{max}} ( - \widehat{\Sigma}^\frac{1}{2} \Lambda \widehat{\Sigma}^\frac{1}{2} ) \Big\}  < \gamma^\star \le \big\| \widehat{\Sigma}^\frac{1}{2} \Lambda \widehat{\Sigma}^\frac{1}{2} \big\|_\ast \Big(\sqrt{\frac{6}{\varepsilon}} \, \mathbbm{1}_{[0, \frac{1}{24}]} + \big(6 + \frac{1}{4\varepsilon}\big) \mathbbm{1}_{(\frac{1}{24},\infty)}\Big)
     \end{array}
    \end{cases}
     \end{align}
    \end{subequations}
    
    \item \textbf{Lipschitz constant:} \label{Prop_Lipschitz_KL}
    The Lipschitz constant of $\dualFunc{\Lambda}$ in Lemma~\ref{lemma_LipschitzContinuity_dLambda} is bounded by
    \begin{equation} \label{Eq_LipschitzConstantKL}
    \Lipschitz \le  \Big(n\sqrt{6 \varepsilon} \, \mathbbm{1}_{[0, \frac{1}{24}]} + n\big(6\varepsilon + \frac{1}{4}\big)\mathbbm{1}_{(\frac{1}{24},\infty)}+1\Big)\|\widehat{\Sigma}\|_{\Fro}.
    \end{equation}    
    \end{enumerate}
\end{proposition}

A similar KL setting as in Proposition~\ref{Prop_Oracle_Lipschitz_KLConst} has been studied in \cite[Proposition 4.1]{Factor_Models_Zorzi}) where a a closed-form solution to the dual of \eqref{Eq_FM_Convex_App} is provided. Proposition~\ref{Prop_Oracle_Lipschitz_KLConst} refines this solution by characterizing the upper and lower bounds~\eqref{Eq_KL_gamma_update} for the optimal dual multiplier. This improvement builds on next lemma that offers a lower bound for $\KL ( \Sigma || \widehat{\Sigma})$.
          
\begin{lemma}[KL-lower bound]\label{Lemma_KL_LowerBound}
The KL divergence $\KL(\Sigma || \widehat{\Sigma})$ in \eqref{Eq_KL} satisfies 
\begin{subequations}
\begin{align}\label{Eq_KL_ineqaulity}
    2 \KL ( \Sigma || \widehat{\Sigma} )  \geq \sum_{i=1}^{n} f\big(\lambda_i(\widehat{\Sigma}^\frac{-1}{2} \Sigma \widehat{\Sigma}^\frac{-1}{2})  \big),
\end{align}
where $ f(\lambda) \coloneqq \frac{1}{3} (\lambda - 1)^2 \mathbbm{1}_{ [0, \frac{3}{2}] } + \big(\frac{1}{3} \lambda - \frac{5}{12}\big) \mathbbm{1}_{(\frac{3}{2} , \infty )}$. For any radius~$\varepsilon \ge 0$, if $\KL\big( \Sigma || \widehat{\Sigma} \big) \le \varepsilon$ we then have 
\begin{align}\label{KL-ball-eps}
\begin{array}{r}
    \big| \lambda_{\rm max}\big(\widehat{\Sigma}^\frac{-1}{2} \Sigma \widehat{\Sigma}^\frac{-1}{2}\big) - 1 \big| \leq \sqrt{6 \varepsilon} \, \mathbbm{1}_{[0, \frac{1}{24}]}
    + \big(6\varepsilon + \frac{1}{4}\big) \mathbbm{1}_{(\frac{1}{24},\infty)}. 
    \end{array} 
\end{align}    
\end{subequations}
\end{lemma}
     
\subsection{Gelbrich distance}

The final case is the Gelbrich distance $\Gel (\Sigma, \widehat{\Sigma})$~\cite{ref:Gelbrich1984}, measuring the distance between two distributions with identical means and covariance matrices $\Sigma$ and $\widehat{\Sigma}$, defined as
\begin{equation}\label{eq:gelbrich}
    \Gel(\Sigma,\widehat{\Sigma}) = \sqrt{\Tr{\Sigma +  \widehat{\Sigma} - 2  \big( \widehat{\Sigma}^\frac{1}{2} \Sigma \widehat{\Sigma}^\frac{1}{2} \big)^\frac{1}{2}}} .
\end{equation} 
We adapt the LMO \eqref{Eq_LMO} under $\Gel (\Sigma, \widehat{\Sigma})$ from \cite{Bridging_Bayesian-Viet}, incorporating a slight generalization that extends to arbitrary matrices (rather than PSD matrices as in \cite{Bridging_Bayesian-Viet}). Additionally, we provide a tight upper bound for the univariate optimization problem and establish the Lipschitz constant of the respective dual function, the two parameters contributing to our algorithm convergence.

\begin{proposition}[Gelbrich oracle \& Lipschitz constant]\label{Prop_Gelbrich} 
Consider LMO \eqref{Eq_LMO} for a given $\widehat{\Sigma} \succeq 0$ where the distance function is the Gelbrich distance $\dist(\Sigma,\widehat{\Sigma})=\Gel(\Sigma, \widehat{\Sigma})$.
\begin{enumerate} [label=(\roman*), leftmargin = 4mm]
    \item\label{Prop_Oracle_Gel}\textbf{Closed-form description:} For any $\Lambda$ and a non-negative scalar $\gamma$, we define 
    \begin{subequations}
    \begin{align} \label{Eq_LMO_Gel}
        \Sigma^\star(\Lambda,\gamma) \coloneqq \gamma^2 (\gamma I +\Lambda)^{-1} \widehat{\Sigma} (\gamma I +\Lambda)^{-1}.
    \end{align}  
    Then, $\LMO(\Lambda)$ in~\eqref{Eq_LMO} equates to $\Sigma^\star(\Lambda,\gamma^\star)$, where $\gamma^\star$ is the unique scalar solution to the concave optimization
    \begin{align}\label{Eq_gamma_Gel}
    \left \{ \begin{array}{cl}
        \max\limits_{\gamma \ge 0} & \gamma \big(\langle  I - \gamma (\gamma I + \Lambda) ^{-1}, \widehat{\Sigma}\rangle  - \varepsilon^2 \big)  \\
        {\rm s.t.} & \lambda_{\rm max}(-\Lambda) < \gamma \le \frac{1}{\varepsilon} \|\Lambda\|_{\rm F}\big(2{\lambda_{\rm max}^{\frac{1}{2}}(\widehat{\Sigma})} + \varepsilon \big) 
    \end{array} \right. 
    \end{align}
    \end{subequations}
    \item\label{Prop_Lipschitz_Gel} \textbf{Lipschitz constant:} The Lipschitz constant of $\dualFunc{\Lambda}$ in Lemma~\ref{lemma_LipschitzContinuity_dLambda} is bounded by         
    \begin{equation} \label{Eq_LipschitzConstantGel}
      \Lipschitz \le \big(2 \lambda_{\max}^{\frac{1}{2}} (\widehat{\Sigma}) + \varepsilon \big) \varepsilon + \|\widehat{\Sigma}\|_{\Fro}
    \end{equation}
    \end{enumerate}
\end{proposition}	

We note that the quasi-closed form description~\eqref{Eq_LMO_Gel} and the lower bound for $\gamma$ in~\eqref{Eq_gamma_Gel} were previously proposed in \cite[Proposition~A.2]{Bridging_Bayesian-Viet}. However, when $\Lambda$ is indefinite, which is an important case in the factor model, the upper bound for $\gamma$ in~\eqref{Eq_gamma_Gel} and the Lipschitz constant~\eqref{Eq_LipschitzConstantGel} are first introduced here. The proof of Proposition~\ref{Prop_Gelbrich} builds on next lemma characterizing the regularity properties of the Gelbrich distance and the respective ball. 

\begin{lemma} [Gelbrich lower bound] \label{lem_GB}
The Gelbrich distance in \eqref{eq:gelbrich} satisfies 
\begin{subequations}
\begin{equation}\label{Gelbrich_LB}
    \Gel (\Sigma , \widehat{\Sigma}) \geq \max\bigg\{ \big\| \lambda^{\frac{1}{2}}(\Sigma) - \lambda^{\frac{1}{2}}(\widehat{\Sigma}) \big\|_2, \, \lambda_{\max}^{-\frac{1}{2}} \big( \Sigma + \widehat{\Sigma} + 2 (\widehat{\Sigma}^\frac{1}{2} \Sigma \widehat{\Sigma}^\frac{1}{2})^{\frac{1}{2}} \big)\| \Sigma - \widehat{\Sigma} \|_{\Fro} \bigg\}.
\end{equation}
Particularly, for any~$\varepsilon \ge 0$, if $\Gel ( \Sigma, \widehat{\Sigma} ) \le \varepsilon$, we have 
\begin{align}\label{Gel_Fer bound}
    \Gel (\Sigma , \widehat{\Sigma}) \geq \frac{\| \Sigma - \widehat{\Sigma} \|_{\Fro}}{2 \lambda_{\max}^{\frac{1}{2}} (\widehat{\Sigma}) + \varepsilon},
\end{align}
\end{subequations}
\end{lemma}

\begin{remark}[Gelbrich strong convexity]\label{rem:strong convexity}
    The inequality~\eqref{Gel_Fer bound} suggests that the Gelbrich distance is strongly convex with respect to the Frobenius norm, uniformly over any compact subset of the PSD cone, a property of interest in view of optimization algorithms. In comparison with the earlier results in the literature~\cite[Theorem\,1]{ref:bhatia2018strong}, the strong convexity relation~\eqref{Gel_Fer bound} does {\em not} depend on the minimum eigenvalue of~$\widehat{\Sigma}$, making it particularly useful when the ball contains low-rank matrices. 
\end{remark}  

\section{Numerical Example} \label{Sec_NumExp}
In our numerical investigation, we evaluate the performance of our algorithm \eqref{Eq_OurFirstOrderAlgorithm} and observe the effect of $\varepsilon$ on covariance matrix estimation accuracy. To this end, we present the numerical results of the algorithm convergence, $\Sigma_\text{True}$ estimation, and the algorithm execution time compared to MOSEK. We initialize our algorithm at a random positive definite matrix $\Lambda_1$ projected to $\SSetProj{1} \cap \SSetProj{2}$ using Dykstra's projection. The algorithm stopping condition is the normalized {\em relative} change as
    \begin{equation}\label{Eq_StoppingCondition}
        \frac{| \langle \Lambda_t , \Sigma_t \rangle - \langle \Lambda_{t-1} , \Sigma_{t-1} \rangle |}{|\langle \Lambda_t , \Sigma_t \rangle|} \leq 10^{-6}.
    \end{equation}
    
    \subsection{Synthetic data generation}\label{Sec_SynDataGen}
        We generate the data in two steps: First, we build the ground-truth covariance matrix $\Sigma_{\textrm{True}} = D_{\textrm{True}} + L_{\textrm{True}}$, where $L_{\textrm{True}} = \Phi_{\textrm{True}} \Phi_{\textrm{True}}^\top$, and $\Phi_{\textrm{True}}$ and $D_{\textrm{True}}$ are created pseudorandomly using \emph{rand} function in MATLAB with \emph{rng(1,'twister')} and \emph{ rng(0)}, respectively. We note that $\Phi_{\textrm{True}}$ and $D_{\textrm{True}}$ are adjusted to be far away from the origin, by setting a lower bound for their minimum element (we set the bound to be $5$, by adding $5$ to \emph{rand} function) so that $\Sigma_\text{True}$ is far from $0$, to exclude the trivial solution $\Sigma = 0$, when $\Lambda \succeq 0$ and the ball includes~$0$. In the next step, we build $N = 15 \hspace{0.5mm} n$ pseudorandom samples of factors $\{\alpha_k\}_{k=1}^N$ and noise $\{\omega_k\}_{k=1}^N$, from zero mean normal distributions with covariance matrices $\Sigma_\alpha= I$ and $D_{\textrm{True}}$, respectively. Then, we build the observed vectors samples $\{ \xi_k = \Phi_\textrm{True} \alpha_k + w_k \}_{k=1}^N$ and compute $\widehat{\Sigma}$.

    \subsection{Convergence}  
    The saddle point problem \eqref{Eq_Our_FM_GeneralDist} is solved using~\eqref{Eq_OurFirstOrderAlgorithm} with $\delta_t= 1/\sqrt{t}$ stepsize, $10^4$ iterations, for the heart disease dataset at \url{www.kaggle.com/datasets/johnsmith88/heart-disease-dataset}, downloaded on April 2026, with $n= 13$ including age, chest pain type, and resting blood pressure to indicate the presence of heart disease, for $\Fro(\Sigma,\widehat{\Sigma})$, $\KL(\Sigma || \widehat{\Sigma})$, and $\Gel (\Sigma,\widehat{\Sigma})$ cases, and for $\varepsilon= \sqrt{10}$,  $\varepsilon= 0.01$, and  $\varepsilon= 0.1$, respectively. Since the true optimal value is unknown, we choose the objective value at iteration $10^4$ as the reference optimal objective value $\langle \Lambda^\star , \Sigma^\star \rangle$ and investigate the algorithm convergence by the normalized convergence error as 
    \begin{align}\label{nor_err}
        e(\Sigma_t):= \frac{|\langle \Lambda_t , \Sigma_t \rangle - \langle \Lambda^\star , \Sigma^\star \rangle|}{|\langle \Lambda^\star , \Sigma^\star \rangle|},
    \end{align}
    We note that to prevent terminating earlier than $10^4$ iterations, stopping condition \eqref{Eq_StoppingCondition} is not activated. Figure~\ref{Fig_conv_and_sweetS}, top row demonstrates that after $200$ iterations, the normalized error~\eqref{nor_err} decreases to around $7.8 \times 10^{-6}$, $0.16$, and $0.02$ of its initial value, for $\Fro(\Sigma,\widehat{\Sigma})$, $\KL(\Sigma || \widehat{\Sigma})$, and $\Gel(\Sigma,\widehat{\Sigma})$ cases, respectively. Hence, the theoretical convergence result of Proposition~\ref{Prop_Alg} is validated. Moreover, based on the red curve in the $\KL$ case, which shows convergence of \cite{Factor_Models_Zorzi} with $\varepsilon= 0.01$, our algorithm outperforms the alternating direction method of multipliers (ADMM) algorithm employed in \cite{Factor_Models_Zorzi}. In this case, the normalized convergence error is defined as $e_t:= \frac{|F(\lambda_t, X_t) - F(\lambda^\star, X^\star)|}{|F(\lambda^\star, X^\star)|}$, with $F$, $\lambda$, and $X$ as in \cite{Factor_Models_Zorzi}. 
    
\subsection{Estimation of the ground-truth~{$\Sigma_\text{True}$}}
\begin{figure*}[t]
	      \centering
	      \includegraphics[width=1\textwidth]{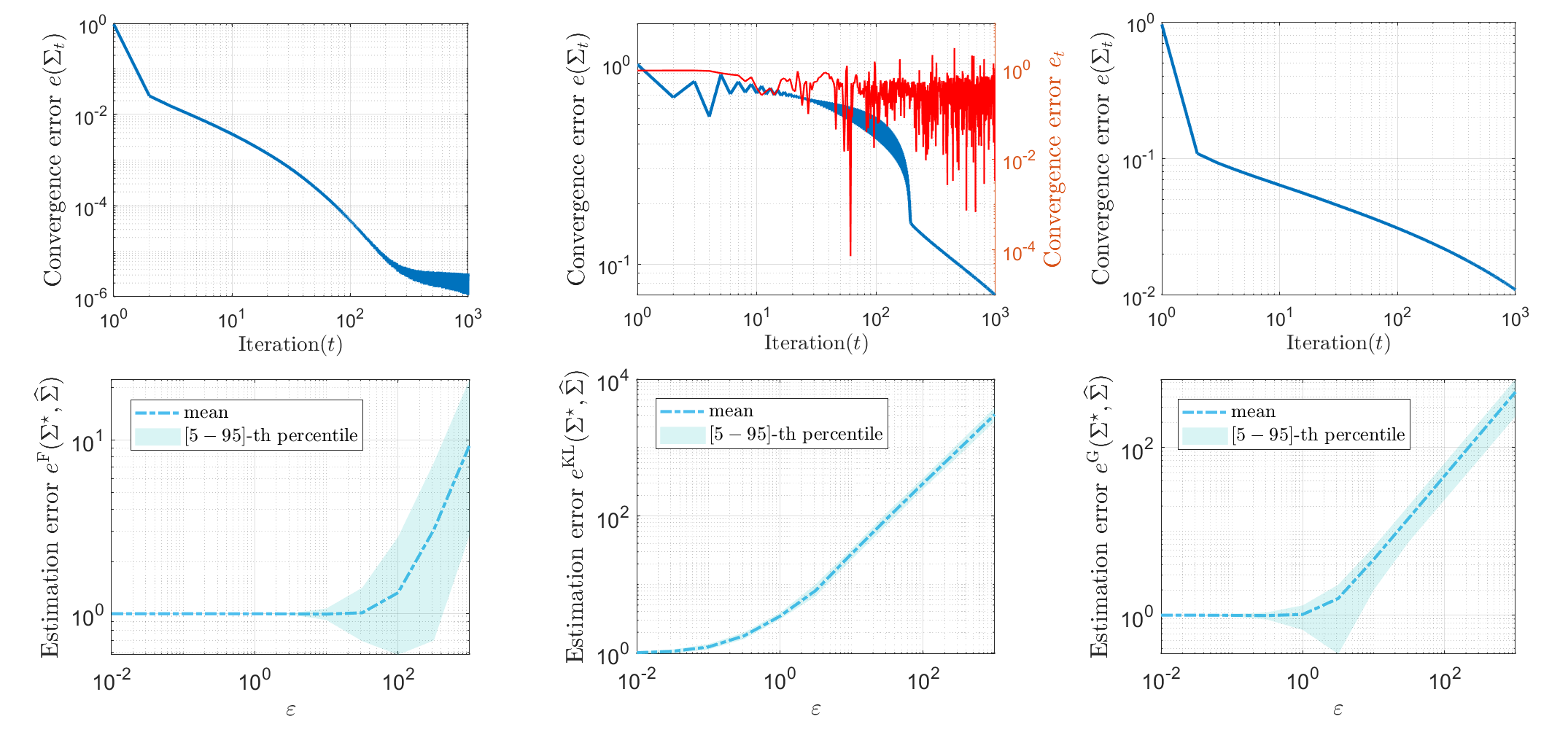}
		    \caption{Analysis of convergence (top row) and estimation of $\Sigma_\text{True}$ (bottom row): Frobenius norm (left column), KL divergence (middle column), and Gelbrich distance (right column).} \label{Fig_conv_and_sweetS}
    \end{figure*}

The effect of the hyperparameter $\varepsilon$ on $\Sigma_{\textrm{True}}$ estimation error is investigated for $N_{exp}= 100$ experiments, for $n= 20$ and $r= 4$. In each experiment, the factor model problem is solved using Algorithm \eqref{Eq_OurFirstOrderAlgorithm} with $\delta_t = 1/\sqrt{t}$ stepsize for $\varepsilon \in \Omega= \{ 0.01(\sqrt{10})^i ~ | ~ i = 0,1,...,10 \}$ and for maximum $t_{\textrm{end}}=10^4$ iterations, while the stopping criterion \eqref{Eq_StoppingCondition} is active. The data is generated through the synthetic data generation procedure in Section \ref{Sec_SynDataGen} with fixed $D_{\textrm{True}}$ and $\Sigma_{\textrm{True}}$ for all experiments. In each experiment, the same samples of observed vectors $\{ \xi_k \}_{k=1}^N$ are used to solve the problem for various~$\varepsilon$, while these samples vary among different experiments. We define the normalized estimation error of $\Sigma_{\textrm{True}}$ as 
\begin{align}\label{Eq_est_err}
    e^{\dist}(\Sigma^\star,\widehat{\Sigma}) \coloneqq \frac{\dist(\Sigma^\star,\Sigma_{True})}{\dist(\widehat{\Sigma},\Sigma_{True})}.    
\end{align}
For $\Fro(\Sigma,\widehat{\Sigma})$ and $\Gel (\Sigma,\widehat{\Sigma})$ cases, an improvement in $\Sigma_\text{True}$ estimation compared to $\widehat{\Sigma}$ is achieved in $61 \%$ and $52 \%$ of the experiments, respectively, as shown as sweet spots in Figure~\ref{Fig_conv_and_sweetS}, bottom row, for $\varepsilon= 100$ and $\varepsilon=\sqrt{10}$. For $\KL(\Sigma||\widehat{\Sigma})$ case, while a sweet spot is not seen, $\Sigma_\text{True}$ estimation is slightly improved in $37 \%$ of the experiments.
\subsection{Execution time} 
\begin{figure*}[t]
      \centering
      \begin{subfigure}{0.32\textwidth}
        \includegraphics[width=\textwidth]{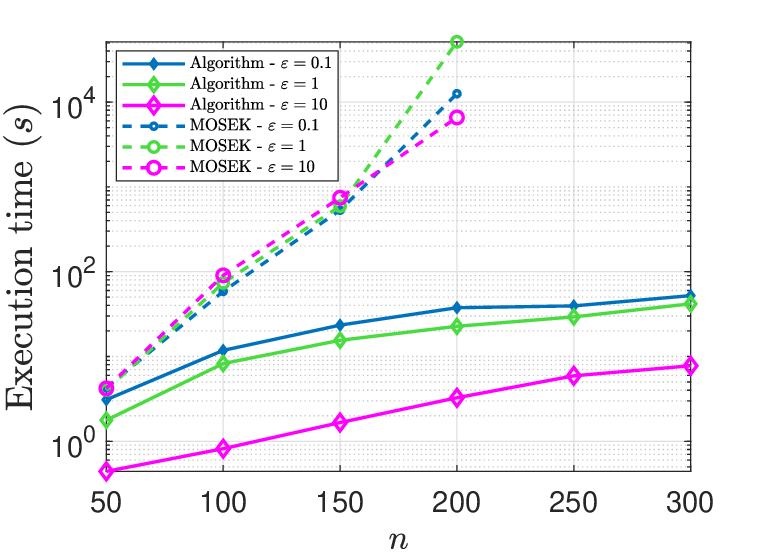}
        \caption{Frobenius norm} \label{Fig_compTimeFro_300}
       \end{subfigure}
       \begin{subfigure}{0.32\textwidth}
    \includegraphics[width=\textwidth]{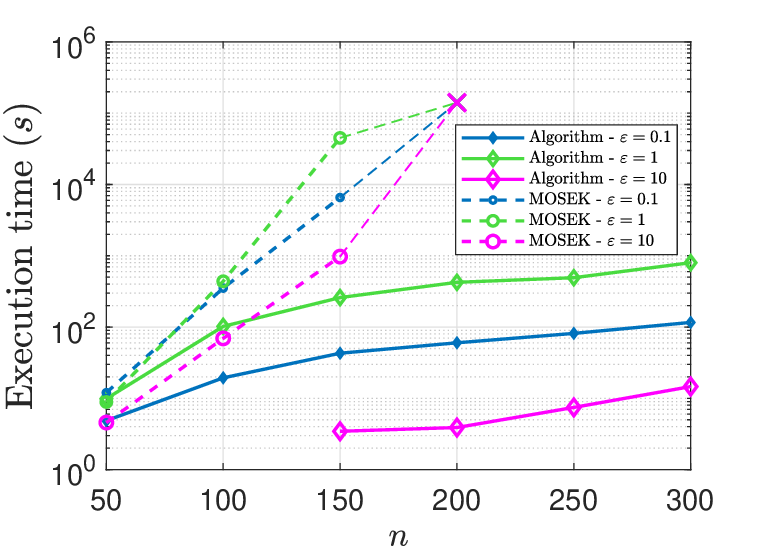}
         \caption{KL divergence\protect\footnotemark} \label{Fig_compTimeKL_300}
        \end{subfigure}
        \begin{subfigure}{0.32\textwidth}
          \includegraphics[width=\textwidth]{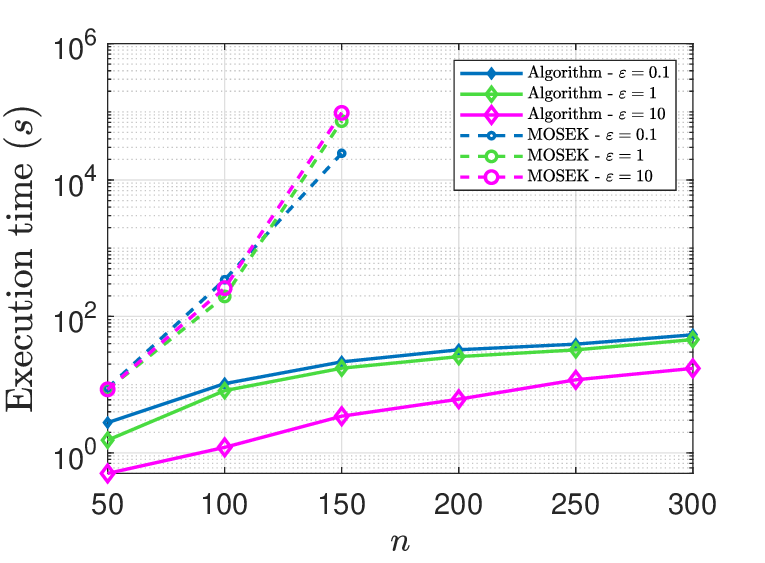}

          \caption{Gelbrich distance} \label{Fig_compTimeGel_300}
        \end{subfigure}
        \caption{Computational time comparison of algorithm \eqref{Eq_OurFirstOrderAlgorithm} and MOSEK.
        \label{Fig_compTime_300}} 
\end{figure*}

For $\Fro(\Sigma,\widehat{\Sigma})$, $\KL(\Sigma||\widehat{\Sigma})$, and~$\Gel(\Sigma,\widehat{\Sigma})$ cases, the execution time of our algorithm with stopping condition~\eqref{Eq_StoppingCondition}, maximum $10^6$ iterations, and $\delta_t = 1/\sqrt{t}$ is compared to MOSEK, for $\varepsilon \in \{0.1, 1, 10\}$ and synthetic datasets generated as in \ref{Sec_SynDataGen}. We solve~\eqref{Eq_FM_Convex_App} and reformulations in Lemmas~\ref{lem:KL-reformulation} and \ref{lem:Gel-reformulation} for $\Fro(\Sigma,\widehat{\Sigma})$, $\KL(\Sigma||\widehat{\Sigma})$, and $\Gel(\Sigma,\widehat{\Sigma})$, respectively, for implementations using MOSEK 9.3 with YALMIP~\cite{YALMIP_Lofberg} in MATLAB R2024a, on a Core(TM) i7-1185G7 CPU with 3.00GHz (1.80 GHz) clock speed and 16GB RAM. According to Figure~\ref{Fig_compTime_300}, our algorithm is clearly more efficient than MOSEK, especially for higher-dimensional data. MOSEK failed to solve the $\Fro(\Sigma,\widehat{\Sigma})$ and $\KL(\Sigma||\widehat{\Sigma})$ cases with $n \geq 250$ and the $\Gel(\Sigma,\widehat{\Sigma})$ case with $n \geq 200$, due to running out of memory. We note that for $\Fro(\Sigma,\widehat{\Sigma})$ and $\Gel (\Sigma,\widehat{\Sigma})$ cases, the execution time for bigger~$\varepsilon$ is less than smaller $\varepsilon$. For $\KL(\Sigma||\widehat{\Sigma})$ case, the effect of increasing $\varepsilon$ on the execution time is not clear.

\section{Conclusion and Future Direction} \label{Sec_Con_FutureDir}
We propose a saddle point reformulation of the factor model problem and a first-order algorithm relying on an LMO. We address the inaccuracy of the empirical covariance matrix~$\widehat{\Sigma}$ by robustification via a family of covariance matrices around it with respect to a generic distance function. We further derive semi-closed form solution for the LMO in three distance functions: Frobenius norm, KL divergence, and Gelbrich distance. For future research, the physical meaning of the factors and factor loading matrix in dynamical systems can be found, which helps to predict system behavior by analyzing these components, or by identifying their variations to investigate stability or detect potential faults. Controllers can then be designed by finding the mapping between control parameters and factor model components.

\footnotetext{~For our algorithm, for $\varepsilon = 10$,  $n=50$ and $n=100$, since the optimal objective value is $0$, the stopping criterion \eqref{Eq_StoppingCondition} is never triggered. Therefore, for these cases we set the iteration limit to $10^4$, yielding optimal objective values of $-3.507$ and $-0.0087$, respectively. These cases where executed in another, similar, hardware. For $n = 200$, MOSEK could not converge after $39$ hours, for all three $\varepsilon$ values.}

\appendix
\section{Technical Proofs}\label{Append_Proofs}
\noindent \emph{Proof of Proposition \ref{Prop_saddle}.~}
Let us consider a new decision variable $\Sigma = L + D$, aligned with~\eqref{Eq_Cov}. The corresponding Lagrange multiplier of this equality constraint is denoted by $\Lambda$, which is symmetric because $\Sigma$, $L$, and $D$ are symmetric. The factor model problem~\eqref{Eq_FM_Convex_App} can then be rewritten as
\begin{align} \label{Eq_FM_Lag} 
\begin{array}{ccl}
    \max\limits_{\Lambda} & \min\limits_{\Sigma, L, D} 
    &  \Tr{L} + \langle \Lambda , \Sigma - L - D \rangle \\
    & \textrm{s.t.} &  L \in \SPlusSet, \quad D \in \DPlusSet, \quad \Sigma \in \Ball{\dist},
  \end{array}
\end{align}
where the strong duality holds thanks to the usual convex-concave property. We note that the Slater's condition is not required due to the linearity of the dualized constraint \cite{ConvexOptimizationTheory_Bertsekas}. The decision variables $L$ and $D$ in~\eqref{Eq_FM_Lag} are separable, and as such, the inner minimization over their respective conic spaces $\SPlusSet$ and $\DPlusSet$ can be computed explicitly. This yields the so-called support functions $\max_{L \in \SPlusSet} \langle \Lambda - I , L \rangle$, and $\max_{D \in \DPlusSet} \langle \Lambda , D \rangle$ in the objective.  It is worth stating that using the dual cone definition, $\max_{L \in \SPlusSet} \langle \Lambda - I , L \rangle = 0 $, if $I - \Lambda \in \SPlusSet^\ast$; otherwise is $+\infty$; and
$\max_{D \in \DPlusSet} \langle \Lambda , D \rangle  = 0$ if $- \Lambda \in \DPlusSet^\ast$; otherwise is $+\infty$; for more general setting, the reader is referred to \cite[Proposition 5.3.9]{ConvexOptimizationTheory_Bertsekas}. Hence, the support functions essentially confine the feasible set of $I - \Lambda$ and $- \Lambda$ to the respective dual cones $\SPlusSet^\ast$ and $\DPlusSet^\ast$, respectively. It then suffices to note that the PSD cone is self-dual (i.e., $ \SPlusSet^\ast=\SPlusSet$) to arrive at the desired program \eqref{Eq_Our_FM_GeneralDist}. \qed

\noindent \emph{Proof of Lemma \ref{lemma_LipschitzContinuity_dLambda}.~}
Using the definition of the dual function in \eqref{Eq_dLambda}, we have
  \begin{align*}
    \dualFunc{\Lambda_1} - \dualFunc{\Lambda_2} & = \min_{\Sigma_1 \in \Ball{\dist}} \langle \Lambda_1 , \Sigma_1 \rangle  - \min_{\Sigma_2 \in \Ball{\dist}} \langle \Lambda_2 , \Sigma_2 \rangle = \min_{\Sigma_1 \in \Ball{\dist}} \max_{\Sigma_2 \in \Ball{\dist}} \langle \Lambda_1 , \Sigma_1 \rangle - \langle \Lambda_2 , \Sigma_2 \rangle \\
     & \leq \max_{\Sigma_2 \in \Ball{\dist}} \langle \Lambda_1 - \Lambda_2 , \Sigma_2 \rangle \leq \Lipschitz \|\Lambda_1 - \Lambda_2\|_{\Fro}, \\
  \end{align*}
where the second equality is derived since for every function $f(x)$, we have $- \min_x f(x) = \max_x -f(x)$, and since $ \langle \Lambda_1, \Sigma_1 \rangle $ is constant when maximizing over $\Sigma_2$. The first inequality holds since the minimum value of the function $\max_{\Sigma_2 \in \Ball{\dist}} \langle \Lambda_1 , \Sigma_1 \rangle - \langle \Lambda_2 , \Sigma_2 \rangle$, over $\Sigma_1$, is less than or equal to its value at any other feasible point including $\Sigma_2$. The last inequality is the direct application of the Cauchy-Schwarz inequality (i.e., $\langle A, B \rangle \leq \|A\|_{\Fro} \|B\|_{\Fro}$ for all $A$, $B$ with appropriate dimensions). \qed

\noindent \emph{Proof of Proposition \ref{Prop_Alg}.~}
The proof follows from the standard projected subgradient method. LMO~\eqref{Eq_LMO} is the subgradient of dual function~\eqref{Eq_dLambda}, as a consequence of Danskin’s theorem~\cite[(A.22), p\,154]{bertsekas1971control}, i.e., $\LMO(\Lambda) = {\frac{\partial \dualFunc{\Lambda}}{ \partial \Lambda}}$. Hence, the first two steps in~\eqref{Eq_OurFirstOrderAlgorithm} effectively implement the projected subgradient ascent applied to~\eqref{Eq_dLambda} while $\Avg{\Lambda}_t$ is the average over all the iterations, i.e., $\Avg{\Lambda}_t = {\frac{1}{N}}\sum_{i=1}^t \Lambda_i$. Leveraging the classical bound of the projected subgradient (e.g., \cite[Proposition\,3.2.4]{bertsekas2015convex}), we have the error bound~\eqref{Alg_bound}. \qed

\noindent \emph{Proof of Proposition \ref{prop_projection}.~}
It suffices to show all requirements of \cite[Theorem~5.3]{wang2024convergence} specified in \cite[Assumption~4.1]{wang2024convergence} are satisfied for $\Pi_{\SSetProj{1} \cap \SSetProj{2}} (\bar{\Lambda})$. Since $\mathbb S_1$ is a polyhedron, it is $C^2$-cone reducible in the sense of \cite[Definition~3.135]{bonnans2013perturbation}. Similarly, $\mathbb S_2$ can be reduced to a (shifted) semidefinite cone. Thus, the first requirement is satisfied.  
    The second and third requirements also hold as ${\text{\normalfont rint}}(\SSetProj{1} \cap \SSetProj{2})$ is nonempty and due to our assumption that $\bar\Lambda - \Lambda^* \in {\text{\normalfont rint}}(\mathcal{N}_{\SSetProj{1} \cap \SSetProj{2}}(\Lambda^*))$. \qed

\noindent \emph{Proof of Proposition \ref{Prop_Oracle_Lipschitz_FroNormConst}.~}
Regarding part $(i)$, by introducing the Lagrange multiplier $\gamma$ for the constraint $\| \Sigma - \widehat{\Sigma} \|_{\Fro}^2 \leq \varepsilon^2$, the dual program of the oracle \eqref{Eq_LMO} is
\begin{align}\label{Eq_LMO_Fro_PrimalDual}
\max_{\gamma \geq 0} \min_{\Sigma \succeq 0} \, \langle \Lambda , \Sigma \rangle + \gamma \big(\| \Sigma - \widehat{\Sigma} \|_{\Fro}^2 - \varepsilon^2 \big). \end{align}
Since the Hessian of the Lagrangian function with respect to~$\Sigma$ is a scaled identity, we first compute the solution of the unconstrained program of the inner minimization (ignoring $\Sigma \succeq 0$), and then project it onto the PSD cone. The first-order optimality condition for the unconstrained inner minimization yields $ \Lambda + 2 \gamma (\Sigma^\ast - \widehat{\Sigma}) = 0$, for every~$(\Lambda, \gamma)$, which describes the unique solution to this problem, yielding the closed-form projected solution $\Sigma^*(\Lambda,\gamma)$ in \eqref{Eq_LMO_Fro_PrimalDual_SigmaS}. Substituting $\Sigma^\star(\Lambda,\gamma)$ in~\eqref{Eq_LMO_Fro_PrimalDual} arrives at~\eqref{Eq_LMO_Fro_gamma} whose solution~$\gamma^\star$  together with the closed-form solution  $\Sigma^\star(\Lambda, \gamma^\star)$ yields the saddle point of~\eqref{Eq_LMO_Fro_PrimalDual}. We note that excluding $\gamma = 0$ is without loss of generality since the objective function of the inner minimization is lower-semicontinuous in $\gamma$ (pointwise minimum of a continuous function in both variables~$(\gamma,\Lambda)$). To show that $\gamma^\star$ is bounded by $\|\Lambda\|_{\Fro}$, we use \cite[Lemma 1]{nedic2009approximate}, which relies on the existence of a Slater point. Considering $\widehat{\Sigma}$ as the Slater point, \cite[Lemma 1]{nedic2009approximate} offers the upper bound 
\begin{equation}\label{Eq_gammaStarUpperBound_Fro}
  \begin{aligned}    
    \gamma^\star & \leq \frac{1}{\varepsilon^2} \Big( \langle \widehat{\Sigma} , \Lambda \rangle - \min_{\Sigma \in \Ball{\Fro}} \langle \Sigma, \Lambda \rangle  \Big) = \max_{\Sigma \in \Ball{\Fro}}{\frac{1}{\varepsilon^2}} \langle \widehat{\Sigma} - \Sigma, \Lambda \rangle \leq \frac{1}{\varepsilon}\|\Lambda\|_{\Fro} 
  \end{aligned}
\end{equation}
where the last inequality is a consequence of lifting the constraint $\Sigma \succeq 0$ and the fact that the Frobenius norm is self-dual. Next, we prove part $(ii)$. Thanks to Lemma 2.2, in particular \eqref{Eq_LipschitzConstant_general}, the Lipschitz constant~$\Lipschitz$ satisfies 
\begin{align*} 
\Lipschitz \le  \max_{\scriptstyle \Sigma \in \Ball{\Fro}} \| \Sigma - \widehat{\Sigma}\|_{\Fro} + \|\widehat{\Sigma}\|_{\Fro} \leq \varepsilon + \|\widehat{\Sigma}\|_{\Fro}
\end{align*}
where the first inequality is the basic triangle inequality of the norm, and the second inequality is the direct consequence of the constraint $\Sigma \in \Ball{\Fro}$. \qed

\noindent \emph{Proof of Proposition \ref{Prop_Oracle_Lipschitz_KLConst}.~}
The proof of both parts follows similar lines as Proposition~\ref{Prop_Oracle_Lipschitz_FroNormConst}. Concerning~(i), dualizing the distance function constraint corresponding to $\KL(\Sigma || \widehat{\Sigma}) \leq \varepsilon$ in the LMO \eqref{Eq_LMO} yields
\begin{equation} \label{Eq_Min_Max_KLConst} 
    \max_{\gamma \geq 0 } \min_{\Sigma \succeq 0} \, \langle \Lambda , \Sigma \rangle  + \frac{\gamma}{2} \Big( - \log \det{\Sigma} + \log \det{\widehat{\Sigma}} + \Tr{\Sigma \widehat{\Sigma}^{-1}} - n \Big) - \gamma \varepsilon\,.
\end{equation}
Due to the term $-\log\det \Sigma$ in the objective of~\eqref{Eq_Min_Max_KLConst}, we can exclude the boundary of the PSD cone for the inner optimizer $\Sigma$ (i.e., $\Sigma \succ 0)$, when $\gamma > 0 $. Similar to proof of Proposition~\ref{Prop_Oracle_Lipschitz_FroNormConst}, we can exclude $\gamma = 0$ since the optimal value of the inner minimization is lower-semicontinuous in $\gamma$. Thus, the first-order optimality condition for the inner minimization yields the unique minimizer~$\Sigma^*(\Lambda,\gamma)$ in \eqref{Eq_KL_SigmaStar}, substituting which in \eqref{Eq_Min_Max_KLConst} yields 
\begin{equation} \label{Eq_Dual_KLConst} 
    \max_{\gamma > 0} \, \langle \Lambda , \Sigma^* (\Lambda,\gamma) \rangle + \frac{\gamma}{2} \Big( - \log \det{\Sigma^* (\Lambda,\gamma)} + \log \det{\widehat{\Sigma}} + \Tr{\Sigma^* (\Lambda,\gamma) \widehat{\Sigma}^{-1}} - n \Big) - \gamma \varepsilon \,.
\end{equation}
Applying the first-order optimality condition to~\eqref{Eq_Dual_KLConst} yields the algebraic equation~\eqref{Eq_KL_gamma_update}. Next, we derive the lower bound for~$\gamma^\star$ solving \eqref{Eq_Dual_KLConst} (or equivalently~\eqref{Eq_KL_gamma_update}). Since $\Sigma^*(\Lambda,\gamma^*) \succ 0$, we have $\widehat{\Sigma}^{-1} + \frac{2}{\gamma^\star} \Lambda \succ 0$, which implies $\gamma^\star >  2 \lambda_{\textrm{max}} ( - \widehat{\Sigma}^\frac{1}{2} \Lambda \widehat{\Sigma}^\frac{1}{2} )$. This bound and the original non-negativity constraint conclude the lower bound. To show the upper bound of~$\gamma^\star$ in \eqref{Eq_KL_gamma_update}, following similar lines as in Proposition~\ref{Prop_Oracle_Lipschitz_FroNormConst}, we have 
\begin{align} \label{bdd_1}
    \gamma^\star &\le \max_{\Sigma \in \Ball{\KL}}{\frac{1}{\varepsilon}} \langle \widehat{\Sigma} - \Sigma, \Lambda \rangle \le 
    \left \{ \begin{array}{rl} 
    \max\limits_{\Sigma \succeq 0} & {\frac{1}{\varepsilon}} \langle \widehat{\Sigma} - \Sigma, \Lambda \rangle \\
   \textrm{s.t.} & \KL \big( \Sigma || \widehat{\Sigma} \big) \le \varepsilon \end{array} \right.
   \le \left \{ \begin{array}{rl} \max\limits_{\Sigma \succeq 0} & {\frac{1}{\varepsilon}} \langle \widehat{\Sigma} - \Sigma, \Lambda \rangle\\
   \textrm{s.t.} & \big| \lambda_{\rm max}\big(\widehat{\Sigma}^\frac{-1}{2} \Sigma \widehat{\Sigma}^\frac{-1}{2}\big) - 1 \big| \\ & \leq \sqrt{6 \varepsilon} \, \mathbbm{1}_{[0, \frac{1}{24}]} + \big(6\varepsilon + \frac{1}{4}\big) \mathbbm{1}_{(\frac{1}{24},\infty)},
   \end{array} \right.
\end{align} 
where the last inequality follows from \eqref{KL-ball-eps}.
The objective function in \eqref{bdd_1} can be upper bounded by
\begin{align*}
   \langle \widehat{\Sigma} - \Sigma, \Lambda \rangle & = \big \langle I - \widehat{\Sigma}^{-\frac{1}{2}} \Sigma \widehat{\Sigma}^{-\frac{1}{2}}, \widehat{\Sigma}^\frac{1}{2} \Lambda \widehat{\Sigma}^\frac{1}{2} \big \rangle \leq \max_{i \leq n} \big | \lambda_i (\widehat{\Sigma}^{-\frac{1}{2}} \Sigma \widehat{\Sigma}^{-\frac{1}{2}}) - 1 \big | \big \| \widehat{\Sigma}^\frac{1}{2} \Lambda \widehat{\Sigma}^\frac{1}{2} \big \|_{\ast},
\end{align*}
Considering the upper bound in~\eqref{bdd_1} and replacing $\lambda_{\rm max}$ yields the desired upper bound for $\gamma^\star$ in~\eqref{Eq_KL_gamma_update}. Regarding $\Lipschitz$ in part~\ref{Prop_Lipschitz_KL}, following bound~\eqref{Eq_LipschitzConstant_general}, we have
\begin{align*} 
    \Lipschitz & = \max_{\Sigma \in \Ball{\KL}} \| \Sigma \|_{\Fro} \le \max_{\Sigma \in \Ball{\KL}} \| \Sigma - \widehat{\Sigma}\|_{\Fro} + \|\widehat{\Sigma}\|_{\Fro} \le \Big(\max_{\Sigma \in \Ball{\KL}} \| \widehat{\Sigma}^{\frac{-1}{2}}\Sigma\widehat{\Sigma}^{\frac{-1}{2}} - I\|_{\Fro} + 1 \Big) \|\widehat{\Sigma}\|_{\Fro} \\
    & \le  \Big(\max_{\Sigma \in \Ball{\KL}}n\lambda_{\rm max}\big(\widehat{\Sigma}^{\frac{-1}{2}}\Sigma\widehat{\Sigma}^{\frac{-1}{2}}-I\big) + 1\Big)\|\widehat{\Sigma}\|_{\Fro} \le \Big(n\sqrt{6 \varepsilon} \, \mathbbm{1}_{[0, \frac{1}{24}]} + n\big(6\varepsilon + \frac{1}{4}\big)\mathbbm{1}_{(\frac{1}{24},\infty)}+1\Big)\|\widehat{\Sigma}\|_{\Fro},
\end{align*}
where the last inequality follows from \eqref{KL-ball-eps}. \qed  

\noindent \emph{Proof of Lemma \ref{Lemma_KL_LowerBound}.~}
To prove \eqref{Eq_KL_ineqaulity}, following \eqref{Eq_KL}, we have
\begin{align*}
    2 \KL \big( \Sigma || \widehat{\Sigma} \big) & = \Tr{\widehat{\Sigma}^{-\frac{1}{2}} \Sigma \widehat{\Sigma}^{-\frac{1}{2}} - I}  - \log \det \big( \widehat{\Sigma}^{-\frac{1}{2}} \Sigma \widehat{\Sigma}^{-\frac{1}{2}} \big) = \sum_{i=1}^{n} \lambda_i(\widehat{\Sigma}^\frac{-1}{2} \Sigma \widehat{\Sigma}^\frac{-1}{2}) - 1 - \log \lambda_i(\widehat{\Sigma}^\frac{-1}{2} \Sigma \widehat{\Sigma}^\frac{-1}{2}) \\
    & \geq \sum_{i=1}^{n} f\big(\lambda_i(\widehat{\Sigma}^\frac{-1}{2} \Sigma \widehat{\Sigma}^\frac{-1}{2})\big),
\end{align*}
where the equalities hold, since the trace and determinant operators are symmetric (i.e., $\Tr{AB}= \Tr{BA}$ and $\det(AB)= \det(BA)$). The inequality follows from the Taylor series expansion with degree 2 of the convex function~$\lambda - 1 - \log(\lambda)$ at $\lambda = 1$ within the interval~$[0,\frac{3}{2}]$, followed by a linear extension of this lower bound for $\lambda > \frac{3}{2}$. To prove \eqref{KL-ball-eps}, we use an inverse function argument for the lower bound function~$f$ in \eqref{Eq_KL_ineqaulity} and apply the inequality~$\max_{i\le n}f(\lambda_i) \le \sum_{i \le n} f(\lambda_i)$. \qed

\noindent \emph{Proof of Proposition~\ref{Prop_Gelbrich}.~}
Regarding~\eqref{Eq_LMO_Gel} and the lower bound of $\gamma$ in~\eqref{Eq_gamma_Gel}, we refer to parts (i) and (ii) in~\cite[Proposition~A.2]{Bridging_Bayesian-Viet}. Concerning the upper bound in~\eqref{Eq_gamma_Gel}, we follow the approach used for $\Fro(\Sigma,\widehat{\Sigma})$ to derive a similar bound to~\eqref{Eq_gammaStarUpperBound_Fro} as
\begin{align*}
   \gamma^\star &\leq \max_{\Sigma \in \Ball{\Gel}}{\frac{1}{\varepsilon^2}} \langle \widehat{\Sigma} - \Sigma, \Lambda \rangle = \left\{ \begin{array}{cl}
       \max\limits_{\Sigma \succeq 0} & \frac{1}{\varepsilon^2} \langle \widehat{\Sigma} - \Sigma, \Lambda \rangle \\
       {\rm s.t.} & \Gel(\Sigma, \widehat{\Sigma}) \leq \varepsilon 
   \end{array} \right. \le \left\{ \begin{array}{cl}
       \max\limits_{\Sigma \succeq 0} & \frac{1}{\varepsilon^2} \| \Lambda \|_{\Fro} \| \widehat{\Sigma} - \Sigma \|_{\Fro} \\
       {\rm s.t.} & \| \Sigma - \widehat{\Sigma} \|_{\Fro}\big(2 \lambda_{\max}^{\frac{1}{2}} (\widehat{\Sigma}) + \varepsilon\big)^{-1} \leq \varepsilon 
   \end{array} \right. \\
   & \le \frac{1}{\varepsilon} \| \Lambda \|_{\Fro} \big(2 \lambda_{\max}^{\frac{1}{2}} (\widehat{\Sigma}) + \varepsilon\big),
\end{align*}
where the second inequality is the application of Cauchy-Schwarz inequality to the objective function and~\eqref{Gel_Fer bound} to $\Gel (\Sigma , \widehat{\Sigma})$ in the constraint. Regarding $\Lipschitz$ in~\eqref{Eq_LipschitzConstantGel}, we apply~\eqref{Gel_Fer bound} to the Lipschitz constant~\eqref{Eq_LipschitzConstant_general} to arrive at 
\begin{equation*} 
\Lipschitz = \max_{\Sigma \in \Ball{\Gel}} \| \Sigma \|_{\Fro} \le  \max_{\Sigma \in \Ball{\Gel}} \| \Sigma - \widehat{\Sigma}\|_{\Fro} + \|\widehat{\Sigma}\|_{\Fro} \leq \big(2 \lambda_{\max}^{\frac{1}{2}} (\widehat{\Sigma}) + \varepsilon \big) \varepsilon + \|\widehat{\Sigma}\|_{\Fro}. 
\end{equation*} \qed

\noindent \emph{Proof of Lemma \ref{lem_GB}.~}
Regarding the first term on the right-hand side of inequality~\eqref{Gelbrich_LB}, note that 
\begin{align} 
     \Gel^2 (\Sigma , \widehat{\Sigma}) & = \Tr{\Sigma +  \widehat{\Sigma} - 2  \left( \widehat{\Sigma}^\frac{1}{2} \Sigma \widehat{\Sigma}^\frac{1}{2} \right)^\frac{1}{2}} = \sum_{i \leq n} \lambda_i(\Sigma) + \lambda_i (\widehat{\Sigma}) \ - 2 \lambda_i^{\frac{1}{2}} (\widehat{\Sigma}^{\frac{1}{2}} \Sigma \widehat{\Sigma}^{\frac{1}{2}}) \notag
     \\ \label{eq:1}
     & \geq \sum_{i \leq n} \lambda_i(\Sigma) + \lambda_i (\widehat{\Sigma}) \ - 2 \lambda_i^{\frac{1}{2}} (\widehat{\Sigma}) \lambda_i^{\frac{1}{2}} (\Sigma) = \sum_{i \leq n} (\lambda_i^{\frac{1}{2}} (\Sigma) - \lambda_i^{\frac{1}{2}} (\widehat{\Sigma}))^2 = \big\| \lambda^{\frac{1}{2}}(\Sigma) - \lambda^{\frac{1}{2}}(\widehat{\Sigma}) \big\|_2^{2} , 
\end{align}
where the inequality follows from the matrix version of Hardy-Littlewood-Polya inequality~\cite[Theorem 3.2]{SomeMatroxRearrangement_Carlen}. As for the second term on the right-hand side of~\eqref{Gelbrich_LB}, 
\begin{align*}
    & \Tr{\big(\Sigma + \widehat{\Sigma} - 2 (\widehat{\Sigma}^\frac{1}{2} \Sigma \widehat{\Sigma}^\frac{1}{2} )^\frac{1}{2} \big) \big( \Sigma + \widehat{\Sigma} + 2 (\widehat{\Sigma}^\frac{1}{2} \Sigma \widehat{\Sigma}^\frac{1}{2} )^\frac{1}{2} \big)} = \Tr{(\Sigma + \widehat{\Sigma})^2 - 4(\widehat{\Sigma}^\frac{1}{2} \Sigma \widehat{\Sigma}^\frac{1}{2})} \\ 
    & = \Tr{\Sigma^2 - 2 \Sigma \widehat{\Sigma} + \widehat{\Sigma}^2} = \Tr{(\Sigma - \widehat{\Sigma})^2}= \| \Sigma - \widehat{\Sigma} \|_{\Fro}^2, 
\end{align*}
where the symmetric property~$\Tr{AB} = \Tr{BA}$ is used to derive the first and second equalities. Using Von Neumann's trace inequality on the above equality gives
$\Gel^2(\Sigma, \widehat{\Sigma}) \lambda_{\max} \big( \widehat{\Sigma} + \Sigma + 2 ( \widehat{\Sigma}^\frac{1}{2} \Sigma \widehat{\Sigma} ^\frac{1}{2} )^\frac{1}{2} \big) \geq \| \Sigma - \widehat{\Sigma} \|_{\Fro}^2$,
which, together with \eqref{eq:1} concludes~\eqref{Gelbrich_LB}. Now observe $\varepsilon \geq \Gel (\Sigma , \widehat{\Sigma}) {\geq} \| \lambda^\frac{1}{2}(\Sigma) - \lambda^\frac{1}{2}(\widehat{\Sigma})\|_2 \geq | \lambda_{\max}^\frac{1}{2}(\Sigma) - \lambda_{\max}^\frac{1}{2}(\widehat{\Sigma})|$
which implies $\lambda_{\max}^\frac{1}{2}(\Sigma) \le \lambda_{\max}^\frac{1}{2}(\widehat{\Sigma}) + \varepsilon$. 
Therefore  
\begin{align*}
  & \lambda_{\max} \big( \Sigma + \widehat{\Sigma} + 2(\widehat{\Sigma}^\frac{1}{2} \Sigma \widehat{\Sigma}^\frac{1}{2} )^\frac{1}{2} \big) \leq \lambda_{\max} \big( \Sigma + \widehat{\Sigma} - 2(\widehat{\Sigma}^\frac{1}{2} \Sigma \widehat{\Sigma}^\frac{1}{2} )^\frac{1}{2} \big) + \lambda_{\max} \big( 4(\widehat{\Sigma}^\frac{1}{2} \Sigma \widehat{\Sigma}^\frac{1}{2} )^\frac{1}{2} \big) \\
   & \leq \Gel^2(\Sigma, \widehat{\Sigma}) + 4 \lambda_{\max}^\frac{1}{2}(\widehat{\Sigma}) \lambda_{\max}^\frac{1}{2}(\Sigma) \leq \varepsilon^2 + 4 \lambda_{\max}^\frac{1}{2}(\widehat{\Sigma}) \big(\lambda_{\max}^\frac{1}{2}(\widehat{\Sigma}) + \varepsilon \big) = \big (2\lambda_{\max}^\frac{1}{2}(\widehat{\Sigma}) + \varepsilon \big)^2 .
\end{align*}
The first two inequalities hold due to sub-additivity and sub-multiplicativity of $\lambda_{\max}$, together with the property~$\lambda_{\max}(A) \leq \Tr{A}$ in the PSD cone.
We apply the above bound in~\eqref{Gelbrich_LB} to arrive at~\eqref{Gel_Fer bound}. \qed

\section{Technical Lemmas to solve \eqref{Eq_FM_Convex_App} using MOSEK}\label{Sec_Apendix}

\begin{lemma}[KL-Factor model reformulation]\label{lem:KL-reformulation}
The factor model problem~\eqref{Eq_FM_Convex_App} with the KL divergence $\dist (\Sigma, \widehat{\Sigma}) = \KL(\Sigma || \widehat{\Sigma})$ is equivalent to the program
\begin{align*}
  \begin{array}{cl}
  \min\limits_{L,\Sigma, Z} & \langle L , I \rangle \\
    {\rm s.t.} & \begin{bmatrix}
     \Sigma & Z \\
     Z^\top & {\rm Diag}(Z)	
   \end{bmatrix} \succeq 0 , ~ L \in \SPlusSet, ~ Z \in \mathbb{L}_n, ~ \Sigma - L \in \DPlusSet \\
    & \log \big(\det{\widehat{\Sigma}}\big) + \Tr{\Sigma \widehat{\Sigma}^{-1}} - n - 2 \varepsilon \leq \sum_i \log \big(Z_{ii}\big),
  \end{array}
\end{align*}
\end{lemma}
For the proof, we refer to~\cite[Section 6.2.3]{SemidefOpt_MOSEK}.

\begin{lemma}[Gelbrich-Factor model reformulation]\label{lem:Gel-reformulation} The factor model problem~\eqref{Eq_FM_Convex_App} with~$\dist (\Sigma, \widehat{\Sigma}) =\Gel (\Sigma, \widehat{\Sigma})$ is equivalent to the LMI program
\begin{align*} 
  \begin{array}{cl}
    \min\limits_{L,\Sigma, C} & \langle L , I \rangle \\
    {\rm s.t.} & \begin{bmatrix}
     \Sigma & C \\
     C^\top & \widehat{\Sigma}	
   \end{bmatrix} \succeq 0	, ~ \Sigma - L \in \DPlusSet , ~ L \in \SPlusSet , ~ \langle \Sigma + \widehat{\Sigma} - 2C , I \rangle \leq \varepsilon^2
  \end{array}
\end{align*}
\end{lemma}
The proof is an application of~\cite[Proposition 2.2]{Bridging_Bayesian-Viet} to~\eqref{Eq_FM_Convex_App}.
            
\bibliographystyle{plain}      
\bibliography{Ref_FactorModel.bib}

\end{document}